\documentclass[11pt]{article}

\usepackage{epsf}
\usepackage{psfrag}
\usepackage{amsmath}
\usepackage{amssymb}
\usepackage{lscape}
\usepackage{latexsym}
\usepackage[usenames]{color}
\usepackage[mathcal]{eucal}
\usepackage{stmaryrd}
\usepackage{textcomp}
\usepackage{setspace} 
\usepackage{graphicx} 
\usepackage{pgfplots} 
\usepackage{epsfig}   
\usepackage{epstopdf} 
\usepackage{appendix}
\usepackage{afterpage}
\usepackage{rotating} 
\usepackage{float} 
\usepackage[misc,geometry]{ifsym} 
\usepackage{comment}

\newtheorem{casestudy}{Case Study}[section]

\newcommand{\ds}{\displaystyle}
\newcommand{\dZ}{{\cal Z \kern -0.7em Z}}
\newcommand{\dC}{{\rm\hbox{C \kern-0.8em\raise0.2ex\hbox{\vrule height5.4pt width0.7pt}}}}
\newcommand{\dQ}{{\rm\hbox{Q \kern-0.85em\raise0.25ex\hbox{\vrule height5.4pt width0.7pt}}}}

\newcommand{\proofbox}{\hspace{\fill}{$\Box$}}

\newcommand\old[1]{}
\newcommand{\beqa}{\begin{eqnarray*}}
\newcommand{\eeqa}{\end{eqnarray*}}

\begin{document}

\title{\bf \Large A Multiobjective Optimization Framework for Irrigation Water Allocation}
\author{
Nahid Sultana\thanks{Department of Business Administration, International Islamic University Chittagong. Email: \texttt{nsultana@iiuc.ac.bd}}
\and
M. M. Rizvi\thanks{Centre for Smart Analytics (CSA), Institute of Innovation, Science and Sustainability, Federation University Australia. Email: \texttt{rizvimath@gmail.com}}
\and
G. M. Wali Ullah\thanks{The University of Chittagong. Email: \texttt{wali@cu.ac.bd} (corresponding author)}
\and
Micah Nehring\thanks{The University of Queensland. Email: \texttt{m.nehring@uq.edu.au} }
}
\maketitle
\noindent\textbf{Abstract}
Sustainable irrigation planning requires balancing economic benefits with environmental flow requirements under increasing climatic and resource constraints. Building on the irrigation optimization framework developed by Ullah and Nehring, this study improves the analytical depth of existing approaches by expanding the feasible decision space and systematically characterizing the full economic--environmental trade-off spectrum. Two single-objective formulations, maximizing net agricultural benefit and minimizing environmental flow deficiency (EFD), are solved to identify boundary solutions that define the limits of the feasible space. These are subsequently integrated into a multiobjective optimization framework using scalarization and evolutionary search techniques to generate a high-resolution Pareto frontier. Numerical experiments on the Muhuri Irrigation Project reveal three key outcomes: (i) a complete scenario view with profits ranging from $0.2 \times 10^{9}$ to $1.497 \times 10^{9}$ and EFD values between 0 and 1200~GL, where 1200~GL represents the theoretical annual maximum under a uniform monthly environmental flow target of 100~GL; (ii) explicit trade-offs demonstrating that higher economic returns are consistently associated with greater ecological shortfalls; and (iii) a computationally efficient approach capable of generating nearly 1000 Pareto-optimal solutions within a short runtime ($\sim$10 seconds), substantially improving solution resolution compared to earlier studies. By transforming point-based optimization into comprehensive trade-off mapping, the proposed framework provides a more informative basis for scenario analysis and decision support in irrigation water allocation, offering a practical extension to existing optimization approaches.


\textbf{AMS} subject classifications. 90B90, 90C29, 90C30, 92B05

\section{Introduction}\label{sec1}
Irrigated agriculture remains one of the largest consumers of freshwater worldwide, particularly in regions where seasonal variability and competing demands place sustained pressure on limited water resources. Although approximately 70\% of the Earth's surface is covered by water, only 2.5\% constitutes fresh water, and less than 1\% of that amount is accessible for industrial, agricultural, and domestic use \cite{Mishra2015}. Of this total, only 0.0002\% is available as surface water runoff, while 0.62\% exists as groundwater, and just 0.001\% is stored in the atmosphere \cite{Matta2010}. This imbalance between water availability and practical usability has made effective irrigation water allocation a critical challenge for ensuring food security, sustaining rural livelihoods, and protecting aquatic ecosystems. In many agricultural systems, irrigation planning must reconcile competing objectives: maximizing crop production while maintaining sufficient river flows to sustain ecological processes. These challenges are intensified by climate change, population growth, and expanding cultivation, which together increase water demand while reducing the reliability of surface water supplies. As a result, irrigation systems increasingly rely on conjunctive use strategies, where surface water is supplemented by groundwater during dry seasons. Under such conditions, allocation decisions are no longer purely operational but become strategic planning problems with long term economic, environmental, and social consequences.
\par
\noindent A significant ecological constraint in this context is the concept of environmental flow (e-flow), refers to the quantity, timing, and quality of water required to sustain the ecological integrity and functioning of a river system. Rather than being a fixed value, environmental flow varies seasonally in accordance with the river’s natural flow regime, supporting key processes such as fish migration, floodplain inundation, wetland regeneration, and groundwater recharge. Excessive water abstraction for agriculture or industry disrupts this natural rhythm, leading to ecosystem degradation, reduced recharge, and long term environmental and social costs. Maintaining minimum river flows is therefore not only an ecological requirement but also an economic and societal necessity. Any shortfall below the required environmental flow level is defined as Environmental Flow Deficiency (EFD), the magnitude of which depends critically on how realistically and flexibly environmental flow targets are specified within planning frameworks.
\par
\noindent These challenges are particularly evident in monsoon dominated agricultural systems. Countries such as Bangladesh experience excessive rainfall during monsoon months alongside acute water shortages in the dry winter season. Irrigation systems in these regions must therefore manage both surplus and scarcity while supporting crop production across highly contrasting hydrological conditions. Allocation decisions directly influence cropping patterns, cultivated area, seasonal yields, and farmer livelihoods, making robust planning tools essential for agricultural sustainability \cite{ Srinivas1994,Ullah2021}. 

    \noindent A major initiative to improve water management in southeastern Bangladesh is the Muhuri Irrigation Project (MIP). The project was developed with a dam regulator system to support dry season irrigation while preventing the intrusion of saline water from the Bay of Bengal into inland areas. However, it has consistently underperformed due to system inefficiencies\cite{Ullah2021}. Effective irrigation systems depend on sound water resource management, and in this context, the MIP was designed to utilize river water for agricultural production through planned and controlled distribution. Over the past several decades, the Bangladesh Water Development Board (BWDB) has implemented numerous water resource development and management programs, among which the Muhuri Irrigation Project represents a targeted effort to enhance agricultural productivity through improved irrigation management.
    Despite its strategic importance at the national level, the MIP has received limited attention in academic research. Consequently, determining the optimal allocation of water among crops and identifying suitable cropping patterns under relevant constraints remains an important and unresolved research challenge.
    \par
\noindent  A key limitation of existing allocation studies for the MIP lies in their restricted representation of decision space. Previous work has evaluated irrigation performance primarily within a narrow range of environmental flow deficits. For example, Ullah and Nehring \cite{Ullah2021} reported a maximum Environmental Flow Deficiency (EFD) of 35 GL under a uniform monthly environmental flow target of 100 GL. Given that EFD is defined as the non negative shortfall between target and actual environmental flow, the theoretical maximum deficit corresponding to zero environmental release would equal the full annual target of 1,200 GL. The reported value therefore represents less than 3\% of the feasible deficit range, indicating that only a limited portion of the potential economic–environmental trade-off space has been explored.
\par
\noindent This narrow focus constrains decision relevance. When irrigation systems are pushed beyond moderate operating conditions such as during prolonged dry seasons or under policy driven environmental constraints planners lack guidance on how economic benefits evolve across the full range of possible environmental shortfalls. As a result, allocation decisions risk understating ecological stress or overestimating economic resilience, limiting their usefulness for real world agricultural planning.
\par
\noindent To address this limitation, this study extends the analytical scope of irrigation allocation by expanding the representation of the economic--environmental trade-off space in the Muhuri Irrigation Project. Unlike previous studies that evaluate system performance within a narrow range of environmental flow conditions, this research systematically explores the full feasible spectrum of Environmental Flow Deficiency (EFD), thereby revealing trade-off relationships that remain unobserved under restricted operating ranges.
\par
\noindent This extension is not merely numerical but analytical in nature. By explicitly characterizing irrigation decisions across a wide continuum of environmental constraints, the study provides a more complete understanding of how economic benefits decline and ecological stress intensifies under progressively constrained water availability. This enables a more realistic and policy-relevant interpretation of irrigation system behavior under both moderate and extreme scarcity conditions.
\par
\noindent Accordingly, this study contributes in two key ways: first, it establishes boundary solutions through single-objective optimization to define the limits of the feasible solution space under consistent physical and operational constraints; and second, it constructs a high-resolution trade-off structure between net economic benefit and environmental flow deficiency, enabling a comprehensive characterization of the multiobjective decision space for irrigation planning.

\par
\noindent The main objective of this research is to enhance the practical applicability and effectiveness of irrigation water allocation models in the Muhuri Irrigation Project. To achieve this, the following specific aims are proposed:
\begin{itemize}
\item[(i)] Develop two single objective mathematical models under a common set of constraints, including pumping balance, minimum area requirements, and total land in the Muhuri Irrigation Project. The objectives are formulated and solved individually: (a) maximization of economic benefit, and (b) minimization of environmental flow deficiency.
\item [(ii)] Compare established multiobjective methods with previous work on the Muhuri Irrigation Project using the same database\cite{Ullah2021}.
\item[(iii)] Formulate single objective mathematical models and obtain optimal solutions using interior point and sequential quadratic programming across different algorithms.
\item[(iv)] Design efficient algorithms, evaluate their performance, and implement a structured solution procedure for the proposed models through computational experiments.
\end{itemize}

\noindent  
\noindent The rest of paper is organized as follows. Section \ref{Insight} discusses the theoretical insights and managerial contributions of this study. Section \ref{Form} reviews previous works and highlights the problems encountered in the Muhuri Irrigation Project.  Section \ref{Moo}  sets out the preliminary concepts used to build the proposed models, especially single objective and multiple objective models. Introduces two new mathematical formulations for constrained optimization problems, including environmental constraints. Section \ref{NuExp} presents numerical experiments and solutions for individual single objective models. Section \ref{MOP}  illustrates the Pareto front obtained from the solutions and compares the results with one another as well as with the previous study. Section \ref{Con} concludes with a summary of findings, limitations, and directions for future research.

\section{Decision and Managerial Insights} \label{Insight}

\noindent This study presents three mathematical optimization models for irrigation water allocation in an irrigation system that relies primarily on surface water while supplementing supply with groundwater during periods of seasonal scarcity. The models are designed to represent alternative planning priorities and to clarify trade-offs between economic performance and environmental flow protection under constrained water availability. Model 1 emphasizes the maximization of net economic benefit, highlighting allocation strategies that favor agricultural production but may increase environmental flow deficit. Model 2 prioritizes the reduction of environmental flow deficiency, supporting ecological sustainability at the expense of lower economic returns. Model 3 integrates both objectives within a multiobjective optimization framework, allowing economic and environmental outcomes to be evaluated simultaneously and supporting balanced allocation decisions.
\par
\noindent A key contribution of this study is the expansion of the decision space beyond the limited operating ranges commonly examined in earlier planning efforts. Previous analysis on MIP largely focused on moderate levels of environmental flow deficit, restricting managerial insight to a narrow set of near optimal outcomes. In contrast, the proposed framework illustrates how net economic benefits evolve across a much wider range of environmental flow deficits within a short evaluation horizon. This extended representation enables farmers and local authorities to assess not only favorable allocation options but also the economic implications of operating under increasingly stressed conditions, which are often encountered during dry seasons.
\par
\noindent The proposed models further enhance decision support by treating both crop land allocation and environmental flow release as endogenous decision variables rather than fixed parameters. This formulation allows the optimization process to determine optimal crop area distribution and environmental flow allocation simultaneously, subject to shared operational constraints such as total cultivable land, pumping balance, and minimum area requirements. As a result, the interaction between agricultural production and ecological protection emerges directly from the model structure rather than being imposed through predefined allocation rules.
\par
\noindent The multiobjective formulation explicitly characterizes trade-offs between net economic benefit and environmental flow deficiency by generating a Pareto set of feasible solutions. Instead of relying on ad hoc weighting schemes, this approach allows decision makers to evaluate multiple compromise strategies and select allocation policies that align with shifting priorities across seasons, hydrological conditions, or policy objectives. By making these trade-offs transparent, the framework supports more defensible and adaptive planning decisions in contexts where economic and environmental goals compete.
\par
\noindent From a managerial perspective, the proposed models function as a practical decision support tool for irrigation planners, water resource managers, and policymakers. The framework enables users to compare allocation strategies across the full operating range, identify zones of diminishing economic returns, and understand the consequences of pushing the system toward higher levels of environmental stress. This reduces reliance on trial and error decision making and supports quicker, more informed responses under water scarcity and institutional pressure.
\par
\noindent The analysis is conducted entirely using secondary data previously published for the Muhuri Irrigation Project, ensuring methodological consistency and direct comparability with earlier studies. All hydrological inputs, crop parameters, and environmental flow requirements are adopted without modification. Accordingly, the contribution of this study lies in improving solution robustness, efficiency, and trade-off representation through alternative optimization formulations rather than in site specific data recalibration or operational redesign.
\par
\noindent The scope of the study is clearly defined. While the models are calibrated for the Muhuri Irrigation Project, the underlying framework is applicable to other regulated irrigation systems that rely primarily on surface water while supplementing supply with groundwater during dry periods. At the same time, the conclusions remain grounded in local hydrological, institutional, and socioeconomic conditions, ensuring relevance without overgeneralization.

\section{Literature Review}\label{Form}
\par
\noindent Efficient irrigation water management is increasingly constrained by the need to simultaneously sustain agricultural productivity and preserve riverine ecosystems. Modern irrigation planning must therefore move beyond profit maximization and explicitly incorporate ecological objectives. A major limitation in existing models is the lack of a direct and quantifiable representation of ecological degradation. Environmental Flow Deficiency (EFD) addresses this gap by defining ecological stress as the non-negative shortfall between target and actual environmental flows. This formulation ensures that only deficits are penalized, avoids overcompensation biases, and enables a transparent integration of ecological thresholds into irrigation allocation models.
\par
\noindent Early irrigation optimization studies primarily focused on economic--hydrological trade-offs without explicitly incorporating ecological constraints. For example, Xevi and Khan (2005) \cite{Xevi2005} developed a multiobjective model maximizing farm income while minimizing irrigation costs and groundwater extraction using a weighted goal programming approach, and Lalehzari et al. (2015) \cite{Lalehzari2015} improved water use efficiency through crop planning and pricing strategies using NSGA-II. However, these approaches implicitly assume that reducing water use automatically preserves environmental flows, which is unrealistic in regulated river systems where flow timing and minimum thresholds are critical. As a result, they fail to capture the non-linear relationship between irrigation withdrawals and downstream ecological degradation.
\par
\noindent Subsequent studies introduced environmental considerations but remained structurally limited. Zeinali et al. (2020) \cite{Zeinali2020} coordinated surface and groundwater use through a WEAP--MODFLOW framework coupled with NSGA-II, while Hoque et al. (2022) \cite{Hoque2022} examined environmental flow requirements in dam-regulated systems. Despite these advances, these models rely on indirect proxies such as groundwater levels rather than explicitly modeling ecological flow deficits, which weakens the enforceability of environmental objectives and limits interpretability in practical decision-making.
\par
\noindent At the operational level, models incorporating field-scale constraints and deficit irrigation strategies improved realism and local efficiency. Wardlaw and Bhaktikul (2004) \cite{Wardlaw2004} integrated soil moisture dynamics and canal capacity, while Sadati et al. (2014) \cite{Sadati2014} demonstrated that deficit irrigation can increase farm income under uncertainty using genetic algorithms. Nguyen et al. (2016) \cite{Nguyen2016} further improved computational efficiency in large-scale systems using ant colony optimization. Unlike system-level approaches such as Ricalde et al. (2022) \cite{Ricalde2022}, these operational models primarily focus on field-scale efficiency rather than basin-wide economic--environmental trade-offs, making them less suitable for allocation problems where environmental flow constraints are binding.
\par
\noindent Crop selection studies reveal further inconsistencies. Mohammadrezapour et al. (2019) \cite{Mohammadrezapour2019} identified wheat as having the highest return-to-water ratio using the Cuckoo Optimization Algorithm (COA), outperforming the Genetic Algorithm (GA) under reduced water allocation, whereas Torabi et al. (2024) \cite{Torabi2024} used Multiobjective Grey Wolf Optimizer (MOGWO) and emphasized integrated management for improving water reliability and reducing aquifer depletion. However, these findings may not generalize across systems because they assume relatively stable water availability and do not fully incorporate dynamic environmental flow requirements. Consequently, recommendations derived from such models may increase ecological stress under water-scarce conditions.
\par
\noindent At broader policy scales, models such as   Bhatia and Rana (2020) \cite{Bhatia2020} used an integrated approach of Linear Programming and Weighted Sum Model (WSM) and Musa (2021) \cite{Musa2021} using Goal Programming operate at aggregated levels. While these approaches provide useful system-wide insights, they do not capture operational-level allocation decisions and therefore cannot fully represent irrigation-level trade-offs between water demand and environmental flow requirements.
\par
\noindent Several studies introduce methodological advances but omit explicit ecological flow representation. Montazar et al. (2010) used Multi-Objective non linear Programming (MOP) \cite{Montazar2010}, Ikudayisi et al. (2018) \cite{Ikudayisi2018} applied Combined Pareto Multi-objective Differential Evolution (CPMDE) and found effective than stochastic multi objective Evolutionary Algorithm named Multiobjective Differential Evolution Algorithm(MDEA),  and Marzban (2021) \cite{Marzban2021} incorporate economic and hydrological considerations but do not include environmental flow deficiency as a binding constraint. This omission limits their ability to evaluate sustainability in ecologically sensitive river systems.
\par
\noindent A significant methodological contribution was made by Lewis and Randall (2017) \cite{Lewis2017}, who linked irrigation decisions with groundwater extraction and crop-specific water requirements which was absent in the model of Xevi and Khan\cite{Xevi2005}. However, their framework relies on conventional multiobjective formulations that do not explicitly identify the extreme boundaries of the feasible solution space, thereby limiting the interpretability of trade-offs between economic and ecological objectives.
\par
\noindent In Bangladesh, Ullah and Nehring (2021) \cite{Ullah2021} applied multiobjective optimization to the Muhuri Irrigation Project (MIP) and followed the model developed by Lewis and Randall\cite{Lewis2017} by using NSGA II. Although their model incorporated environmental flow targets, it explored only a limited portion of the feasible trade-off space and required substantial computational effort. This limitation is particularly critical in the MIP, where seasonal variability, tidal influences, and saline intrusion create highly constrained allocation conditions. Many existing optimization frameworks may not fully generalize to such systems due to assumptions of stable inflows or weak coupling between irrigation withdrawals and ecological thresholds.
\par
\noindent Recent studies using advanced computational techniques introduce additional limitations. Ricalde et al. (2022) \cite{Ricalde2022} applied NSGA-III for climate adaptation planning, while Li et al. (2011) \cite{Li2011} developed uncertainty-based optimization frameworks. These approaches provide valuable system-level insights but remain computationally intensive.
\par
Similarly, Chang et al. (2024) \cite{Chang2024} and Mirzaie et al. (2024) \cite{Mirzaie2024} proposed advanced fuzzy multiobjective optimization models incorporating physical processes and environmental factors. While these approaches enhance realism, they may not systematically characterize the complete Pareto-optimal frontier. Moreover, evolutionary algorithms typically produce approximate solution sets without explicitly guaranteeing convergence to exact boundary solutions, which may affect consistency and reproducibility in policy-relevant applications.
\par
\noindent Existing studies can therefore be broadly categorized into three groups: (i) economic-hydrological optimization models that neglect explicit ecological constraints, (ii) environmentally aware models that lack tractable and interpretable solution structures, and (iii) computational frameworks that do not fully characterize the economic-ecological trade-off space. These approaches generally do not simultaneously minimize environmental flow deficiency, identify extreme boundary solutions, and ensure computationally efficient and reproducible optimization for large-scale irrigation systems such as the Muhuri Irrigation Project.
\par
\noindent This study departs from conventional multiobjective formulations by independently solving two single-objective optimization problems, maximizing economic benefit and minimizing environmental flow deficiency under identical constraints. This extreme-point analysis explicitly defines the bounds of the feasible solution space, enabling a more explicit interpretation of trade-offs compared to heuristic Pareto approximations.
\par
\noindent Furthermore, deterministic optimization methods, including interior-point and sequential quadratic programming, are employed. Unlike evolutionary algorithms, these methods offer strong convergence properties, improved computational efficiency, and reproducible solutions under well-defined conditions, making them suitable for policy-sensitive irrigation systems.
\par
\noindent By combining boundary solution identification with efficient deterministic optimization, this study provides insights that are not readily observable through many conventional approaches. Specifically, it reveals how irrigation allocation decisions interact with binding environmental flow constraints under highly constrained hydrological conditions.
\par
\noindent Finally, the proposed framework integrates economic returns, crop-specific water requirements, seasonal variability, and explicit environmental flow constraints within a unified structure tailored to the Muhuri Irrigation Project, enabling direct applicability to real-world irrigation management and bridging the gap between theoretical optimization and operational decision-making.
\section{Problem Statement}\label{paramodel}
\par
Despite these objectives and its socioeconomic importance, the performance of the MIP has fallen short of initial expectations. During the dry winter months, surface water availability from the Muhuri River declines sharply, while rainfall contributes only about 2\% of the annual total \cite{Ullah2021}. Although farmers tend to prefer Rabi crops because of their lower water requirements compared to Kharif crops, the limited availability of surface water during this period restricts the effective use of irrigation. As a coping strategy, many farmers increasingly depend on groundwater extraction, which raises production costs and places additional pressure on already stressed aquifers.
\par
\noindent Institutional and infrastructural limitations further exacerbate these challenges. Inadequate water distribution systems, poor drainage conditions, and unplanned cropping intensities reduce the effectiveness of available irrigation water. Moreover, the centralized operation of the canal network, combined with limited control over on farm water delivery, often results in inequitable access to water, particularly disadvantaging farmers located at the tail ends of canals. Competition among users intensifies during dry periods, increasing the likelihood of inefficient and uneven water use. From an ecological perspective, reduced river inflows caused by irrigation withdrawals contribute to Environmental Flow Deficiency (EFD), threatening aquatic ecosystems and local fisheries that depend on sufficient freshwater availability. At the same time, excessive reliance on groundwater undermines long term sustainability by accelerating aquifer depletion and increasing energy consumption.
\par
\noindent These interconnected socioeconomic, institutional, and environmental challenges highlight the limitations of existing operational practices within the Muhuri Irrigation Project. They underscore the need for an updated and integrated planning framework that explicitly accounts for both economic objectives and ecological constraints. A robust decision support system capable of guiding irrigation water allocation under conditions of seasonal scarcity is therefore essential to improve equity, provide a comprehensive representation of Environmental Flow Deficiency, support informed decision making by local authorities, enhance agricultural productivity, and ensure the long term sustainability of the Muhuri Irrigation Project.
\section{Formulations of the Mathematical Models and Solution Approach}\label{Moo}
An irrigation planning authority must allocate land and water resources efficiently across multiple crops and time periods, ensuring economic viability and environmental sustainability. Surface water is delivered through canals, and additional water can be pumped. The system must meet staple crop area and production demands while minimizing environmental flow deficiency (EFD).\\

\noindent 
Decision variables:
\begin{table} [H]
\begin{tabular}{l l}
$X_c$& Crop Area for each crop c in hectare $ha$\\
$Env.Flow_m$& Environmental Flow in month m in gigaliter $GL$\\
\end{tabular}
\end{table}
\newpage
\noindent Input parameters:
\begin{table} [H]
\begin{tabular}{l l}
$P_c$& Price of total crop production in $AUD$ per hectare (AUD/ha)\\
$C_w$& Cost of Surface Water conveyance (AUD/GL)\\
$K_{cm}$  & Crop coefficient in the different growth stages of crops c in month m\\
$ET_m$& Evapotranspiration($GL/ha$) in each month m\\
$R_m$  & Rainfall($GL/ha$) in each month m \\
$Allocation_m$& Surface Water($GL$) in each month m\\
$P_m$  & Ground or Pumping Water($GL$) in each month m\\
$C_p$& Cost of pumping Groundwater(AUD/GL)\\
$Vcost_c$& Variable cost for each crop c per hectare (AUD/ha)\\
$Inflow_m$& Total runoff or Available River water ($GL$)in each month  \\
$Tar.Env.Flow_m$ & The minimum required Environmental flow($GL$) in each month \\
$T_{\text{pump}}$ & Total Pumping Water  $(GL)$\\
$T_{\text{area}}$ & Total Area $(ha)$\\
$min_{area}$ & Minimum Area $(ha)$
\end{tabular}
\end{table}

\noindent Decision variables dependent parameters:
\begin{table} [H]
\begin{center}
\begin{tabular}{l l}
$Allocation_m$& Allocated Surface water, where $Allocation_m= {Inflow}_m -  {Env.Flow}_m$\\
$P_m$& Pumped Ground Water, where $ P_m = \sum_{c} (K_{cm} \cdot ET_m - R_m  ) X_c - Allocation_m$\\
\end{tabular}
\end{center}
\end{table}
\noindent The objective functions can be formulated as follows:

\noindent The first objective, $f_1$:
\begin{multline}\label{revobj}
    NB_{\text{max}} =  \sum_{c} P_{c} X_c - C_w \cdot  \sum_{m} \left(\sum_{c} (K_{cm} \cdot ET_m -R_m) X_c -  P_m \right)\\
    - C_p \cdot \sum_{m} P_m - \sum_{c} Vcost_c \cdot X_c
\end{multline}
The second objective, $f_2$:
\begin{equation}\label{envobj}
    EFD_{\text{min}} = \sum_{m} \max \left( \text{Tar.Env.Flow}_m - \text{Env.Flow}_m,\ 0 \right)
\end{equation}

\noindent {\bf Model 1: When Net Benefit is maximized}
\begin{equation}\label{eqmathemodel1}
\begin{array}{rl} 
\ds \max_{(X_c, Env.Flow_m) \in X} & \ f_1(X_c, Env.Flow_m)\\
\mbox{subject to the constraints} \\
& \ds X_c, Env.Flow_m \geq 0.
\end{array}
\end{equation}
where the feasible set $X$ is defined as 
\[\hspace*{-5mm} X := \left
\{ (X_c, Env.Flow_m) \in \mathbb{R}^{c} \times \mathbb{R}^{m} \mid  \sum_{m} P_m \leq T_{pump}, \;\;min_{area} \leq X_c,\;\; \sum_{c} X_c \leq T_{area}
\right \}.\]
\noindent The model seeks to maximize net benefit, with $X_c$ and $Env.Flow_m$ as the decision variables and parameters, including 
$P_c$, $C_w$, $K_{c,m}$, $ET_m$, $R_m$, $Allocation_m$, $P_m$, $C_p$, and $Vcost_c$. 
In this setup, maximum net benefit is achieved by sacrificing environmental performance, illustrating the trade-off between economic gain and ecological sustainability. 

\noindent {\bf Model 2: When Environmental Flow Deficiency is minimized}
\begin{equation}\label{eqmathemodel2}
\begin{array}{rl} 
\ds \min_{(X_c, Env.Flow_m) \in X} & \ f_2(X_c, Env.Flow_m)\\
\mbox{subject to the constraints} \\
& \ds X_c, Env.Flow_m  \geq 0.
\end{array}
\end{equation}

\noindent Model~2 is formulated as a single-objective optimization problem in which minimizing
Environmental Flow Deficiency (EFD) is the sole objective. The target environmental flow
$\text{Tar.Env.Flow}_m$ represents a minimum ecological requirement rather than an upper
bound. Accordingly, the decision variable $\text{Env.Flow}_m$ may exceed this target when
sufficient water is available or when environmental protection is prioritized. Any release
above the target does not contribute to EFD, which is defined only as the shortfall below
$\text{Tar.Env.Flow}_m$. Once the minimum environmental flow requirement is satisfied, EFD
becomes zero for all feasible solutions.

\noindent For this reason, Model~2 is not intended to generate a realistic cropping plan in isolation.
Instead, it serves as a boundary single-objective model that identifies the ecologically optimal
extreme of the feasible solution space. In this formulation, crop allocation is governed by
feasibility constraints such as minimum crop area requirements, total land availability,
and pumping capacity rather than by economic incentives. The resulting solution is used as
a reference point for evaluating trade-offs between economic and environmental objectives
in the multiobjective formulation of Model~3.

\noindent {\bf Model 3: Multiobjective - Net Benefit versus Environmental Flow Deficiency }
\begin{equation}\label{eqmathemodel3}
\begin{array}{rl} 
\ds \min_{(X_c, Env.Flow_m) \in X}  & \ds \ \left\{-f_1, f_2\right\}\\
\mbox{subject to the constraints} \\
& \ds X_c, Env.Flow_m \geq 0.
\end{array}
\end{equation}

\noindent
This multiobjective programming model simultaneously minimizes two conflicting goals: 
the negative Net Benefit function $f_1(X_c,Env.Flow_m)$ and the environmental flow deficiency function $f_2(X_c,Env.Flow_m)$. 
The constraints ensure feasible pumping limits, enforce both minimum and maximum land use requirements, 
and guarantee that crop area and environmental flow allocations remain non-negative. 
In essence, the formulation captures the trade-off between maximizing agricultural returns 
and preserving ecological water requirements.

\noindent We develop code within the MATLAB environment and utilise MATLAB solvers to address the models. A comprehensive solution procedure and algorithm for solving Models 1–3 are detailed in Appendices \ref{solprocedure} and \ref{mopalg}.

\section{Numerical Experiments,  Results and Sensitivity Analysis}\label{NuExp}

\noindent 
This section presents a numerical study demonstrating that optimal irrigation water allocation can be achieved under a large number of operational constraints during dry year conditions. To account for the stochastic nature of the optimization solvers and their sensitivity to initial conditions, multiple randomized trials are performed. The solution procedure described in \eqref{solprocedure} is applied to solve the single objective formulations, where Model~1~\eqref{eqmathemodel1}, maximizes net economic benefit, and Model~2 ~\eqref{eqmathemodel2}, minimizes environmental flow deficiency. The single objective problems are solved using MATLAB nonlinear optimization solvers, including `fmincon' implemented with the interior point method and the sequential quadratic programming (SQP) algorithm. For the multiobjective formulation, MATLAB’s genetic algorithm and Pareto search approach are employed, together with the scalarization method described in \cite{Burachik2014}, to approximate the Pareto front and explore trade-offs between economic and environmental objectives. These solvers are evaluated based on their ability to handle large scale optimization problems and the challenges associated with approximating the Pareto front. All optimization experiments are conducted with randomly generated initial decision vectors and repeated multiple times to ensure robustness. For each solver and model combination, computational time and solution quality are recorded and summarized using representative statistics. All numerical computations are carried out on an HP Evo laptop equipped with $16$~GB RAM and a Core~$i7$ processor operating at $4.6$~GHz.

\par
\begin{casestudy}\label{ex1} 
The numerical example is based on the framework developed by Ullah and Nehring \cite{Ullah2021},
who applied an optimization model to the Muhuri Irrigation Project using dry-year data
representative of low-flow, water-stressed conditions in Bangladesh. In the present study,
dry-year scenarios are interpreted as operational worst-case conditions for irrigation
planning, rather than as indicators of arid or climatically comparable environments.
Accordingly, the modeling assumptions adopted here ensure methodological consistency
and comparability with prior studies, while remaining aligned with the hydrological and
institutional setting of the Muhuri Irrigation Project.
\par
\noindent The analysis considers a regional water authority responsible for managing the ten most
demand-intensive crops over a twelve-month planning horizon. The total cultivable area
is 23,076 hectares, with a maximum annual groundwater pumping capacity of 50~GL and a
minimum crop area requirement of 1,000 hectares per crop. The unit costs of pumped and
surface water are set at 100,000 and 26,000~AUD per GL, respectively. The use of AUD,
rather than Bangladeshi Taka (BDT) or USD, is maintained to ensure direct comparability
with prior work. The target environmental flow is specified as 100~GL per month.(A fixed monthly environmental flow target is used to maintain consistency with the
benchmark framework and to enable transparent comparison across alternative
optimization formulations. Although environmental flow requirements in monsoon-
dominated river systems are inherently seasonal, introducing seasonally varying targets
would add hydrological complexity beyond the scope of this methodological study.
Future work may incorporate seasonally adaptive environmental flow requirements.)
\par
\noindent Key parameter values, including pumping limits and environmental flow targets, are adopted
directly from existing studies \cite{Ullah2021, WaliUllah2020} to enable consistent benchmarking of alternative
optimization formulations under identical operating conditions. These parameters are not
intended as site-specific policy prescriptions, but rather as representative constraints for
evaluating solution behavior and trade-off characteristics. All hydrological inputs, including river inflow, rainfall, and evapotranspiration, are
adopted directly from the published datasets used in Ullah and Nehring~\cite{Ullah2021},
which were compiled from official project records of the Muhuri Irrigation Project. These
data are treated as exogenous inputs to ensure methodological consistency and to focus
the analysis on optimization behavior rather than hydrological estimation. Within this
framework, allocation outcomes are examined both with and without enforcing the
environmental flow constraint to identify the upper bound of attainable economic benefit
and to explore trade-offs between profit maximization and environmental flow protection.

\end{casestudy}
 
	\begin{table} [H]
 		\caption{\small{\textit{Optimal Allocation of Crops and Environmental Flow for Model 1 \eqref{eqmathemodel1}. }}}
		\footnotesize
		\vskip 1.5em
 		\centering
 		\begin{tabular}{|c| c| c| c| c| c| c| c| c| c| c| c|c|}
 			\hline
 			crops &   Aus &  Aman  & Boro & Wheat& Potato & Oil &Pulses & Suger & Winter & Summer & &  \\ 
             &   rice &  rice  & rice & &  & seeds & & cane & vegi & vegi &  &  \\ [0.5ex]
 			\hline
$X_c$&1000 &1000 & 1000 &1000&2076&1000&1000&5000&5000&5000& &\\ [.5ex]\hline
Month&Jan &Feb & Mar &April&May&June&July&Aug&Sep&Oct&Nov & Dec\\ [.5ex]\hline
E.Flow&0 &0 & 0  &0&0&0&0&0&0&0&0 &0\\ [.5ex]\hline
EFD&100 &100 & 100  &100&100&100&100&100&100&100&100 &100\\ [.5ex]\hline
                \end{tabular}
		\label{table:ex1a}
	\end{table}
 \par

\noindent We solved Model 1 \eqref{eqmathemodel1} and obtained the crop and environmental flow allocations shown in Table~\ref{table:ex1a}. In this case, the maximum profit achieved is $ \$ 1.4968 \times 10^9$. However, since the second objective of maintaining environmental flows was fully relaxed, no water was allocated to environmental requirements, and no pumping was required. As a result, the environmental flow deficit reaches its maximum possible value of 1200 GL
\par
\noindent Now we consider that the environment flow meets the target flow requirement of the environment in Model~2~\eqref{eqmathemodel2}, and Model~ 3~\eqref{eqmathemodel3}, that is,

\begin{equation}\label{tarcon}
    \ds Tar.Env.Flow_m \leq Env. Flow_m,
\end{equation}

	\begin{table} [H]
 		\caption{\small{\textit{Optimal crop and environmental flow allocation for Model 2 \eqref{eqmathemodel2} under constraint \eqref{tarcon}. }}}
		\footnotesize
		\vskip 1.5em
 		\centering
 		\begin{tabular}{|c| c| c| c| c| c| c| c| c| c| c| c|c|}
 			\hline
 			crops &   Aus &  Aman  & Boro & Wheat& Potato & Oil &Pulses & Suger & Winter & Summer & &  \\ 
             &   rice &  rice  & rice & &  & seeds & & cane & vegi & vegi &  &  \\ [0.5ex]
 			\hline
$X_c$&1000 &1000 & 1000 &1000&2076&1000&1000&5000&5000&5000& &\\ [.5ex]\hline
Month&Jan &Feb & Mar &April&May&June&July&Aug&Sep&Oct&Nov & Dec\\ [.5ex]\hline
E.Flow&100 &100 & 100  &100&100&100&100&100&100&100&100 &100\\ [.5ex]\hline
EFD&0 &0 & 0  &0&0&0&0&0&0&0&0 &0\\ [.5ex]\hline

                \end{tabular}
		\label{table:ex1b}
	\end{table}
 \par

\noindent We solved Model~2~\eqref{eqmathemodel2} using constraint \eqref{tarcon} and obtained the crop and environmental flow allocations presented in Table~\ref{table:ex1b}. In this scenario, the maximum profit is $ \$1.4080 \times 10^9$. Here, the environmental flow deficiency is $0$ GL, since constraint \eqref{tarcon} does not allow for any shortfall.  To satisfy this requirement, additional water must be supplied through pumping when surface water is insufficient.
\par
\noindent The crop area allocations reported in Tables~1 and~2 exhibit relatively tidy values, with
several crops assigned areas close to the minimum allowable level. This outcome is a
direct consequence of the modeling structure rather than an artefact of the solution
process. In particular, the presence of minimum crop area constraints, together with
upper bounds on individual crop areas and the absence of any objective term promoting
crop diversification, leads the optimization algorithms to select boundary solutions.
Once feasibility conditions are satisfied, continuous nonlinear solvers such as
`fmincon' typically converge to corner points of the feasible region. As a result,
the reported allocations should be interpreted as constraint-driven extreme solutions
that define the bounds of the feasible space, rather than as detailed operational cropping
plans.

	\begin{table} [H]
 		\caption{\small{\textit{Optimal crop and environmental flow allocation for Model~3  \eqref{eqmathemodel3}. }}}
		\footnotesize
		\vskip 1.5em
 		\centering
 		\begin{tabular}{|c| c| c| c| c| c| c| c| c| c| c| c|c|}
 			\hline
 			crops &   Aus &  Aman  & Boro & Wheat& Potato & Oil &Pulses & Suger & Winter & Summer & &  \\ 
             &   rice &  rice  & rice & &  & seeds & & cane & vegi & vegi &  &  \\ [0.5ex]
 			\hline
$X_c$&1000 &1000 & 1000 &1000&1000&1000&1000&1000&1000&1000& &\\ [.5ex]\hline
Month&Jan &Feb & Mar &April&May&June&July&Aug&Sep&Oct&Nov & Dec\\ [.5ex]\hline
E.Flow&200 &200 & 200  &200&200&200&200&200&200&200&200 &200\\ [.5ex]\hline
EFD&0 &0 & 0  &0&0&0&0&0&0&0&0 &0\\ [.5ex]\hline
                \end{tabular}
		\label{table:ex1c}
	\end{table}
 \par
\noindent Next we minimize Model~3~\eqref{eqmathemodel3} with the same parameter settings as in Case Study~\ref{ex1}, considering both the inclusion and exclusion of constraint~\eqref{tarcon}. The optimal crop and environmental flow allocations, shown in Table~3, are identical in both situations. In this case, the total environmental flow deficiency is zero, and the profit obtained is $ \$0.2742 \times 10^9$. In Table~3, the environmental flow exceeds the target value of 100~GL per month and is uniformly
allocated at 200~GL, reflecting surplus environmental releases under the multiobjective formulation
rather than a change in the target environmental flow requirement.
 Note: The computational time taken to solve all three cases is approximately 8–10 seconds.
\par
\noindent \textbf{Observations from the Results:}
The three scenarios deliver different outcomes because each goal controls how water is shared between crops and the river:
\begin{itemize}
    \item Profit only: When the model focuses only on profit, nearly all available water is directed to the most valuable crops. This generates the highest income of $ \$ 1.4968 \times 10^9$, but the river receives too little water, creating a total environmental flow shortage of $1200$ GL.
    \item Profit with environmental target: When the model must also meet the river’s target flow, it reserves enough water for the environment. With less surface water left for irrigation, more pumping is needed, which is expensive. As a result, profit drops moderately to $\$ 1.4080 \times 10^9$, but the environmental shortage is reduced to zero.
    \item Minimum environmental shortage: When protecting the environment is the top priority, the model sends a steady 200 GL each month to the river before allocating water to crops. Because of this and limits on pumping, not all land can be irrigated, and profit falls sharply to $\$0.2742 \times 10^9$. The river’s flow target, however, is fully met.
\end{itemize}
\noindent Overall, these results highlight a clear trade-off: directing more water to crops increases income but harms ecological health, while stronger environmental protection lowers profit.

\section{Multiobjective Optimization Problems}\label{MOP}

 \noindent Multiobjective optimization deals with finding the best balance between two or more objectives that often compete with each other. Because these goals conflict, there is rarely a single solution that is best for every objective at once. Instead, the aim is to identify a group of equally valid trade-off solutions called the Pareto front, where improving one objective would worsen at least one other. A common way to address such problems is scalarization, which reformulates a multiobjective problem as a single-objective one so that standard optimization techniques can be applied. In this work, we use the weighted-constraint method \cite{Burachik2014}. Here, the objective functions are expressed with adjustable weights: one weighted objective is optimized, while the remaining objectives are treated as constraints whose bounds depend on the chosen weights. By varying these weights evenly, we can approximate the complete Pareto optimal set. We chose the weighted constraint approach for its straightforward implementation and its ability to generate accurate Pareto solutions (more scalarization approaches in \cite{BurKayRiz2017}, \cite{PascolettiSerafini1984}). The detailed procedure for the two objective case is described below.

\noindent \textbf{\em The weighted-constraint approach ($P_w^k$)\em\,:} For a positive weights we define

\[\hspace*{20mm} W := \left
\{ w \in \mathbb{R}^{2} \mid w_1, w_2 > 0,\sum_{i=1}^{2} w_i=1
\right \}.\]
Let $w \in W$, $f_1$ in \eqref{revobj}, and $f_2$ in \eqref{envobj}, and we therefore define two scalar problems as
\begin{equation} \label{nchep1}
\mbox{($P_{w}^1$)}\ \left\{\begin{array}{rl} \ds\min_{(X_c, Env.Flow_m) \in  X} & \
\ w_1f_1(X_c, Env.Flow_m),
\\[4mm]
\mbox{subject to} & \ \ w_2f_2(X_c, Env.Flow_m) \leq w_1f_1(X_c, Env.Flow_m),\\ & \ \ \ds X_c, Env.Flow_m  \geq 0.
\end{array}
\right.
\end{equation}
and
\begin{equation} \label{nchep2}
\mbox{($P_{w}^2$)}\ \left\{\begin{array}{rl} \ds\min_{(X_c, Env.Flow_m)\in  X} & \
\ w_{2}f_{2}(X_c, Env.Flow_m),
\\[4mm]
\mbox{subject to} & \ \ w_1f_1(X_c, Env.Flow_m) \leq w_2f_2(X_c, Env.Flow_m), \\ & \ \ \ds X_c, Env.Flow_m  \geq 0.
\end{array}
\right.
\end{equation}

\noindent In Section \ref{NuExp}, net profit and environmental deficiency are optimized separately, and the results highlight the trade-off between these two objectives. In Case Study \ref{ex2}, we solve the two objective optimisation problem defined in  \eqref{eqmathemodel3}.
\color{black}
\begin{casestudy}\label{ex2} 
In this analysis, we use the input parameters from Case Study \ref{ex1} to solve the multiobjective optimization problem defined in Model 3. Because it is a large scale problem, different approaches are applied, including the 'gamultiobj' and 'paretosearch' solvers in the MATLAB environment, while the weighted-constraint method is employed to approximate the Pareto front.
\end{casestudy}

\begin{figure}[H]
\begin{minipage}{80mm}
\begin{center}
\includegraphics[width=79mm]{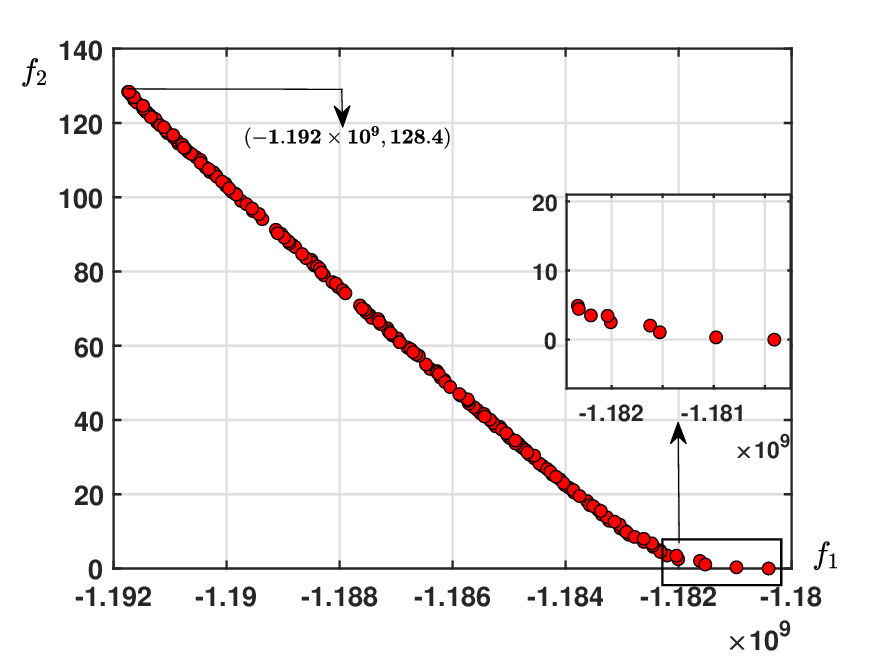} \\
{\scriptsize (a) Pareto front approximated with 'gamultiobj'.}
\end{center}
\end{minipage}
\begin{minipage}{80mm}
\begin{center}
\includegraphics[width=80mm]{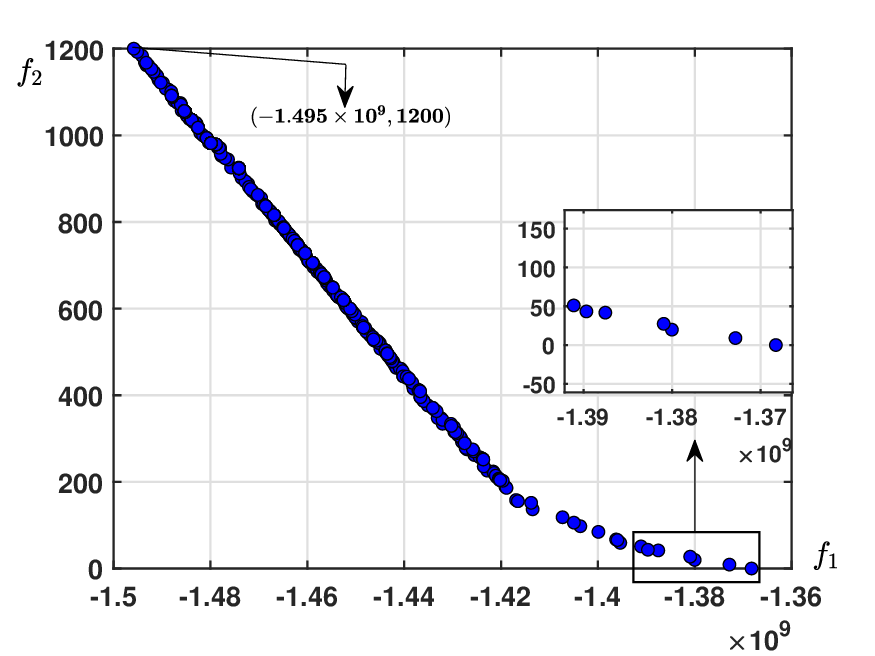} \\
{\scriptsize (b) Pareto front approximated with 'paretosearch'.}
\end{center}
\end{minipage}
\\[2mm]
\begin{minipage}{80mm}
\begin{center}
\includegraphics[width=80mm]{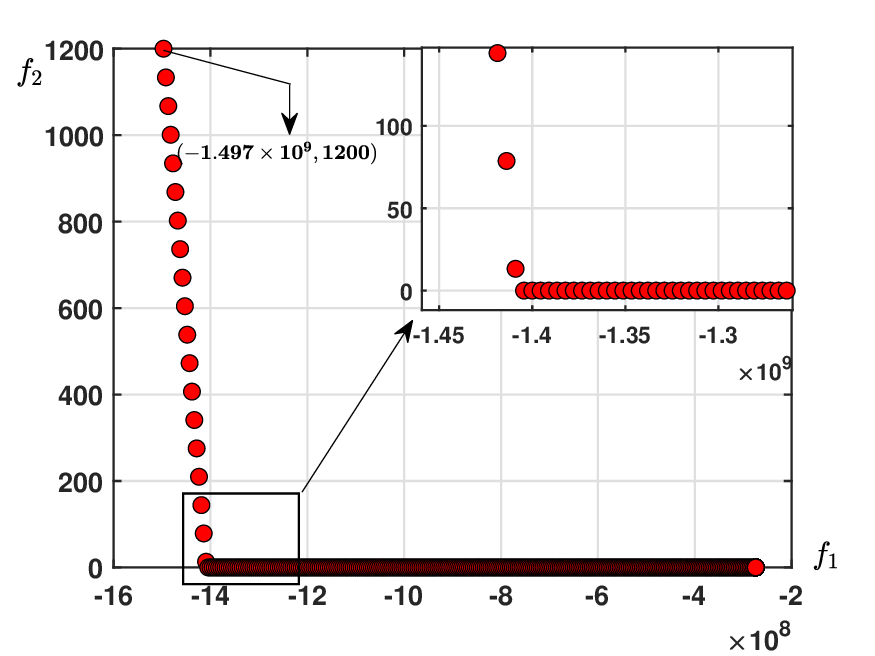} \\
\scriptsize (c) Pareto front approximated with ($P_w^k$).
\end{center}
\end{minipage}

\caption{{Pareto fronts generated by the 'gamultiobj', 'paretosearch', and 'weighted-constraint' ($P_w^k$) methods. }}
\label{multiexample}
\end{figure}

\noindent We solve the problem described in \eqref{eqmathemodel3} using the `gamultiobj' solver, increasing the population size up to $500$. The lower bounds are set as the minimum area for $X_c$ and zero for environmental flow, while the upper bounds are set at $5,000$ ha for $X_c$, (The per-crop upper bound is implemented as it ensures realistic crop distribution while preventing extreme monocropping and preserving comparability with the benchmark framework) and $300$ GL for environmental flow. The Pareto front obtained with this algorithm is shown in Figure~\ref{multiexample}(a), and the computation time along with the total number of Pareto points is reported in Table \ref{tablemop}.

\noindent We now solve the problem defined in \eqref{eqmathemodel3} using the `paretosearch' solver, with the population size increased to 500. The same lower and upper bounds are considered for $X_c$ and environmental flow values. The Pareto front approximated by this algorithm is shown in Figure~\ref{multiexample}(b), while the computation time and the total number of Pareto points are reported in Table \ref{tablemop}.

\noindent Finally, we use the weighted constraint scalarization approach. This method requires a sequence of steps detailed in Appendix \ref{mopalg}. With weight normalisation as defined in $W$, the procedure involves solving subproblems \eqref{nchep1} and \eqref{nchep2} and comparing their solutions, discarding any non-dominated points. To solve the subproblems, we employ the fmincon solver with the `Interior Point', `SQP', and `SQP Legacy' algorithms, using the same lower and upper bounds as in the previous two cases. Following the algorithm outlined in Appendix \ref{mopalg}, we then solve the problem defined in \eqref{eqmathemodel3} to approximate the Pareto points, which are illustrated in Figure~\ref{multiexample}(c). A numerical summary is provided in Table \ref{tablemop}.

\begin{table}[ht]
\small
\caption{Numerical performance of 'gamultiobj', 'paretosearch' and $P_w^k$.}  
\centering 
\begin{tabular}{||c| c| c| c| c| c| c||} 
\hline \hline 
Methods &   Number of  & Population & Number of   & Time          & Pareto    & Non-Pareto  \\  
       &   grid points  & size   & subproblems & in seconds    & points    &   points \\ [0.5ex] 
\hline 
gamultiobj & --   & 500    & --       & 15.72       & 345        & Nil     \\ [.5ex]\hline 
paretosearch            & --       & 500  &  --    & 16.86       & 500        & Nil     \\ [.5ex]\hline
$P_w^k$            & 500   &  --  & 1000        & 10.12       & 989        & 11     \\
\hline
\hline 
\end{tabular}
\label{tablemop} 
\end{table}

%

\noindent \textbf{Comparison Analysis:} We observe that the `gamultiobj' algorithm generates $345$ Pareto points in $15.72$ seconds under the same parameter settings. From Figure~\ref{multiexample}(a), the highest environmental deficiency level is around $130$ GL when profit is $\$ 1.91 \times 10^9$, and the deficiency drops to zero when profit is between ($\$ 1.180 \;\mbox{and}\; 1.181) \times 10^9\). The `paretosearch' algorithm takes slightly more time but produces more Pareto points $500$ in total. In generating the Pareto front, what matters more than the number of points is how much of the Pareto front each algorithm can cover. For paretosearch, the highest deficiency is $1200$ GL at a profit of $\$ 1.49 \times 10^9$, with deficiency reaching zero when profit is around $\$ 1.37 \times 10^9$ (see, Figure~\ref{multiexample}(b)). Using the weighted-constraint scalarization approach, we find that the highest deficiency and corresponding profit are the same as 'paretosearch'. However, the deficiency starts to appear when the profit is $\$ 1.4 \times 10^9$ and decreases to $\$ 0.2 \times 10^9$ (see, Figure~\ref{multiexample}(c)). 

\noindent To evaluate the true Pareto front, we refer to Case Study \ref{ex1} results, where individual objectives $f_1$ and $f_2$ are optimised. The results show that the highest deficiency of $1200$ GL occurs at a profit of $\$ 1.49 \times 10^9$, consistent with Figures~\ref{multiexample}(b) and (c) obtained for `paretosearch' and the `scalarization' approach, respectively. Moreover, maximising $f_1$ under target environmental flow constraints \eqref{tarcon} in Case Study 5.1 results in zero deficiency when profit is $\$ 1.4 \times 10^9$, which matches only the scalarization results in Figures~\ref{multiexample}(c). This indicates that `paretosearch' and `gamultiobj' fail to approximate the lower part of the Pareto front. Therefore, considering computation time, the number of generated points, and the accuracy of the true Pareto front, the scalarization approach performs better than the other methods.\\

\noindent \textbf{Comparison analysis with the previous study\cite{Ullah2021}:} The findings of the present study highlight several critical improvements over the earlier work of Ullah and Nehring \cite{Ullah2021}. A key point of divergence lies in the treatment of Environmental Flow Deficiency (EFD). Ullah and Nehring reported a maximum EFD of only 35 GL under a fixed monthly environmental flow target of $100$ GL. Given that EFD is defined as the shortfall between the prescribed environmental flow and the actual allocation, the theoretical maximum deficit for a uniform 100 GL monthly target is $1200$ GL per year (see, obtained results in Table \ref{table:ex1a}). This indicates that their model captured only a narrow portion of the feasible solution space. Consequently, the ecological shortfalls are understated, and the trade-offs between agricultural benefits and ecological requirements are not fully explored. In contrast, our scalarization based approach succeeds in covering the entire feasible frontier, generating Pareto solutions that extend up to the theoretical maximum EFD. This broader coverage, as illustrated in Figure~\ref{multiexample}(c) (Pareto fronts generated by the Weighted-Constraint approach, $P_w^k$), provides a far more comprehensive representation of the trade-off surface compared to \cite[Figure~1 in Page 10]{Ullah2021}, where the solution set remains confined to a very limited region.For cross method comparison, relative performance is discussed based on the shape and extent of the attainment surfaces rather than absolute magnitudes alone.
\noindent The two studies also differ markedly in terms of computational performance and solution diversity. Ullah and Nehring, using the NSGA-II evolutionary algorithm, reported only $34$ Pareto optimal solutions, obtained at the expense of approximately $706.59$ minutes of computational time. In sharp contrast, the scalarization technique produces $989$ Pareto solutions in just $10.12$ seconds. This demonstrates not only substantial computational efficiency but also a significant improvement in the diversity of solutions available to decision makers. The larger and more representative solution set provided by our approach ensures a more reliable understanding of trade-offs and expands the range of strategies that can be considered in irrigation planning.

\noindent Together, these results establish that the scalarization framework 
achieves superior performance in both accuracy and efficiency. By covering the full spectrum of environmental deficits 
and delivering a wide set of optimal solutions with minimal computational burden, our study offers a more rigorous and practical foundation for decision making in the management of the Muhuri Irrigation Project.

\section{Conclusion}\label{Con}
\par
\noindent This study advances irrigation water allocation analysis by extending the analytical depth of existing optimization frameworks applied to the Muhuri Irrigation Project. Rather than introducing a new model, the contribution lies in expanding the feasible decision space and systematically characterizing the full economic--environmental trade-off between net economic benefit and Environmental Flow Deficiency (EFD).
\par
\noindent The results demonstrate that earlier analyses captured only a limited portion of the feasible trade-off region, whereas the proposed approach identifies the complete EFD spectrum ranging from 0 to 1200~GL, corresponding to the theoretical annual maximum under a uniform monthly environmental flow target. By solving single-objective formulations, the study explicitly defines boundary solutions representing the extreme limits of economic maximization and environmental protection. These boundary points are subsequently used to construct a high-resolution Pareto frontier, enabling continuous and transparent interpretation of trade-offs across the entire feasible range.
\par
\noindent The multiobjective analysis further shows that increasing economic returns are consistently associated with higher environmental deficits, while stricter environmental protection leads to measurable reductions in economic benefit. In addition, the scalarization-based approach generates a dense set of Pareto-optimal solutions (approximately 1000 points) within a short computational time, substantially improving the resolution and completeness of trade-off representation compared to earlier studies. This enhanced resolution provides decision-makers with a broader and more flexible set of allocation strategies under varying environmental constraints.
\par
\noindent From a practical perspective, the framework transforms irrigation optimization from point-based solution analysis into comprehensive trade-off mapping, offering a more informative basis for scenario evaluation and policy-oriented decision support in water allocation under scarcity.
\par
\noindent The study is subject to several limitations. The analysis is conducted within a deterministic framework using fixed input data and dry-year conditions, and therefore does not capture uncertainties associated with climate variability, crop yield fluctuations, or dynamic environmental flow requirements. In addition, socio-economic and institutional factors influencing water allocation decisions are not explicitly incorporated.
\par
\noindent Future research can extend this work by integrating stochastic or fuzzy representations of hydroclimatic uncertainty, incorporating seasonally varying environmental flow targets, and modeling additional hydrological processes such as soil moisture dynamics, reservoir storage, and groundwater interactions. Further development of integrated decision-support systems that combine optimization with socio-economic and policy modules would enhance the applicability of the framework for real-world irrigation management.

\noindent \textbf{\large {Acknowledgment:}} 
\noindent We acknowledge ChatGPT (5) AI tools, which are used to rephrase, edit, and polish the author's written text for spelling, grammar, or general style. 

\section{\small Funding and/or Conflicts of interests/Competing interests}
 The authors declare that no funds, grants, or other support were received during the preparation of this manuscript. The authors have no relevant financial or non-financial interests to disclose.



\section{Appendix}\label{appen}

\subsection{Solution procedure}\label{solprocedure}
\noindent The steps below are used to solve Models $1$ and $2$, as described in Section \ref{NuExp}.
 \begin{description}
\item[Step $\mathbf{1}$] { \textbf{(Input)}} \\ 
The dataset used in our analysis, taken from Ullah and Nehring \cite{Ullah2021}, provides hydrological and agricultural data for the Muhuri Irrigation Project. The parameters considered in the model include $P_c$, $C_w$, $K_{c,m}$, $ET_m$, $R_m$, $Inflow_m$, $C_p$, $Vcost_c$, $Tar.Env.Flow_m$, $T_{pump}$, $T_{area}$, and $min_{area}$ (refer to the parameter description in Section \ref{Form}). 
 \item[Step  $\mathbf{2}$] {\textbf{(Define Problem)}}\\
Define the constraint problem stated in \eqref{eqmathemodel1} to maximize revenue profit and in \eqref{eqmathemodel2} to minimize environmental deficiency. Use the `fmincon` solver with default options and random initial points, and use four algorithms: Active Set, Interior-Point, SQP, and SQP-Legacy. Set lower and upper bounds for both $X_c$ and $Env.Flow_m$.
 \item[Step  $\mathbf{3}$] \textbf{(Solve Problem)}\\
Solve the problems $\ds \max f_1$ and $\ds \min f_2$, subject to the constraints defined in \eqref{eqmathemodel1} and \eqref{eqmathemodel2}. Record the solutions $\left(\bar{X^1}_{c}, \bar{Env.Flow^1}_{m}\right)$. 
 \item [Step  $\mathbf{4}$] \textbf{(Algorithm Validation)}\\
Various algorithms, including Active Set, Interior Point, SQP, and SQP Legacy, are implemented in the MATLAB environment to solve Models 1 and 2 and validate the optimal solutions. The performance of these algorithms is assessed by comparing the results, evaluating computational time, and examining the accuracy of the solution approximations.
\end{description}

\subsection{Algorithm}\label{mopalg}
\noindent The following algorithm is used in Case Study-\ref{ex2} to generate the Pareto fronts for objectives $f_1$ and $f_2$, respectively.

\noindent \textbf{Algorithm 1} \label{algo1}
\begin{description}
	\item[Step $\mathbf{1}$] { \textbf{(Input)}} \\ Set all parameters that are stated in the Case Study-\ref{ex1}.
	\item[Step  $\mathbf{2}$] {\textbf{(Determine the individual minima)}} Solve $\ds \min f_1$ and $\ds \min f_2$, subject to the constraints given in Model \eqref{eqmathemodel3} for the parameters specified in Case Study-\ref{ex1}, and obtain solutions $\left(\bar{X}_{c}, \bar{Env.Flow}_{m}\right)$, respectively.	
 
	\item[Step $\mathbf{3}$] { \textbf{(Generate weighted parameters)}} \\
In this step, we generate the weight parameters $w_i$. Using the weight generation method described in \cite[Step 3 of Algorithm 3]{BurKayRiz2022}, the weights are constructed based on the boundary solution obtained in Step~2.
	\item[Step $\mathbf{4}$] Choose $w=(w_1,1-w_1)$, which generated from Step 3.
	\begin{description}
		\item[(a)] Solve Subproblem \eqref{nchep1} and Subproblem \eqref{nchep2}. Obtain $\bar{x}_1:=\left(\bar{X^1}_{c}, \bar{Env.Flow^1}_{m}\right)$ and $\bar{x}_2:=\left(\bar{X^2}_{c}, \bar{Env.Flow^2}_{m}\right)$, respectively. 
		\item[(b)] Determine Pareto points : 
				\begin{description}
			\item[(i)] If $\bar{x}_1=\bar{x}_2$, then set $\bar{x}=\bar{x}_1$ (a Pareto point)\\
			and record the point.
			\item [(ii)] If $\bar{x}_1=\bar{x}_2$  does not hold, then any dominated point is discarded by comparing these $2$ solutions. \\
			Record the non-dominated Pareto points.	
		\end{description}
	\end{description}
	\item[Step $\mathbf{5}$] (Output)\\
	All recorded points are Pareto points of Case Study-\ref{ex2}.\\
\end{description}

\noindent \textbf{Author Contributions}: 
All authors contributed equally to this work and share the same responsibility for its content.

 
\end{document}